\documentclass[10pt]{article}
\usepackage[utf8]{inputenc} 
\usepackage{geometry} 
\usepackage{graphicx} 
\usepackage[parfill]{parskip} 
\usepackage{booktabs} 
\usepackage{array} 
\usepackage{paralist} 
\usepackage{verbatim} 
\usepackage{subfig} 
\usepackage{amsmath}
\usepackage{amssymb}
\usepackage{mathtools}
\usepackage{txfonts}
\usepackage{enumerate}
\usepackage{float}
\makeatletter
\renewcommand\section{\@startsection{section}{1}{\z@}%
                                  {-3.5ex \@plus -1ex \@minus -.2ex}%
                                  {2.3ex \@plus.2ex}%
                                  {\normalfont\large\bfseries}}
\makeatother

\begin{document}

\title{A Note on Quartic Equations with only Trivial Solutions }
\author{Felix Sidokhine}
\maketitle

\section{Introduction}

This notice is the result of a series of experiments that were inspired by some of the theorems presented in Mordell’s text ``Diophantine Equations". In his treatise of quartic diophantine equations with only trivial solutions, he showed that the triviality of solution of $x^4 + y^4 = 2z^2$ with constraint $(x, y) = 1$ the equation $x^4 + 6x^2y^2 + y^4 = z^2$ subject to the same constraint also has only trivial solution over integers \cite{Mordell1969}. While his text gives this result as an isolated phenomenon, we have decided to research it in a deeper context by studying the equivalence (in terms of solvability) between other quartic equations using a tool that we shall call a ``resolvent". For quartic equations, the resolvent is a system of equations of second-degree which is algebraically obtained from the original diophantine equation. From the collected experimental data, the resolvent should generally have the form: 

\begin{equation}
\begin{cases}
mx^2 + ny^2 = kx'^2 + ly'^2 \\
xy = x'y' \\ 
(x,y)=(x',y')=1
\end{cases}
\end{equation}

So far this approach has generated some interesting results describing the relationship between various fourth-degree diophantine equations over integers and some promising data for the case of factorial extensions of integers (e.g. $\mathbb{Z}[\sqrt{2}]$ and also $\mathbb{Z}[i]$).
\\
\\
In this notice, which carries an expository character, we decided to show our findings by presenting the reader with a very easy to follow and transparent example of how the proposed technique works. In particular we shall show the equivalence and triviality of solutions of the equations $x^4 + 4y^4 = z^2$ and $x^4 + 6x^2y^2 + y^4 = z^2$ subject to the constraint $(x, y) = 1$, and then present the data we obtained from our entire experiment.

\section{Reduction to the Resolvent}

First, we shall show that both of the equations described previously reduce to the same resolvent.

\subsection{$x^4 + 4y^4 = z^2$}

Given the equation $x^4 + 4y^4 = z^2$ we want to build the resolvent that contains the solvability information. Let $(x_0,y_0,z_0)$ be a solution of $x^4 + 4y^4 = z^2$. Without loss of generality we can assume that $(x_0,2y_0)=1$, otherwise if they had a g.c.d the equation simply returns to its original form or the 4 moves from $y$ to $x$ (this can be verified by direct computation).
\\
\\
Therefore we have $x_0^4 + 4y_0^4 = z_0^2$, or $(x_0^2)^2 + (2y_0^2)^2 = z_0^2$. Then $(x_0^2,2y_0^2,z_0)$ is a primitive pythagorean triplet and we can claim that there exist $u_0, v_0$, relatively prime such that:
\begin{equation}\label{eq1}
\begin{dcases}
x_0^2 = u_0^2 - v_0^2 \\
2y_0^2 = 2u_0v_0 \\
z_0 = u_0^2 + v_0^2
\end{dcases}
\end{equation}

Since $\mathbb{Z}$ is a UFD, (\ref{eq1}) implies that $u_0 = s_0^2$ and $v_0 = t_0^2$ for some relatively prime $s_0,t_0$. Substituting in the $x_0$ equality in (\ref{eq1}), we obtain that : $x_0^2 = s_0^4 - t_0^4$
\\
\\
Using the difference of squares identity  we have: 
\begin{equation}\label{eq2}
x_0^2 = (s_0^2 - t_0^2)(s_0^2 + t_0^2), 
\end{equation}
where $t_0 \equiv 0 \mod 2$. Due to $(s_0,t_0)=1$, we also have that $(s_0^2-t_0^2,s_0^2 + t_0^2)=1$.
\\
Once more by using the fact that $\mathbb{Z}$ is a UFD combined with (\ref{eq2}), there are such $\alpha_0, \beta_0$ such that:
\begin{equation}\label{eq3}
\begin{cases}
\alpha_0^2 = s_0^2 + t_0^2 \\
\beta_0^2 = s_0^2 - t_0^2 \\
(s_0,t_0) =1
\end{cases}
\end{equation}

Once more, by Pythagorean triplets and (\ref{eq3}), we obtain the following two statements:

\begin{equation}\label{eq4}
\alpha_0^2 = s_0^2 + t_0^2 \rightarrow
\begin{cases}
s_0 =  \lambda_0^2 - \gamma_0^2 \\
t_0 = 2 \lambda_0 \gamma_0
\end{cases}
\end{equation}

\begin{equation}\label{eq5}
\beta_0^2 = s_0^2 - t_0^2  \rightarrow
\begin{cases}
s_0 =  \lambda_0'^2+ \gamma_0'^2 \\
t_0 = 2 \lambda_0' \gamma_0'
\end{cases}
\end{equation}

Where $(\lambda_0, \gamma_0) = (\lambda_0', \gamma_0') = 1$. Combining (\ref{eq4}) and (\ref{eq5}), we obtain the following result:
\begin{equation}
\begin{cases}
\lambda_0^2 - \gamma_0^2 = \lambda_0'^2+ \gamma_0'^2 \\
\lambda_0 \gamma_0 = \lambda_0' \gamma_0' \\
(\lambda_0, \gamma_0) = (\lambda_0', \gamma_0') = 1
\end{cases}
\end{equation}
Which implies that the system of equations:
\begin{equation}\label{resolvent1}
\begin{cases}
X^2 - Y^2 = X'^2+ Y'^2 \\
X Y = X' Y' \\
(X,Y) = (X', Y') = 1
\end{cases}
\end{equation}
 
is solvable over $\mathbb{Z}$ and most important for this notice, that if $x^4 + 4y^4 = z^2$ was solvable, then so was this system of equations. The system (\ref{resolvent1}) is called the resolvent and can is denoted $\text{Res}(x^4 + 4y^4 = z^2)$.
\\
\\
To show the converse, perform the following substitution:
\begin{eqnarray}\nonumber
T_0 = \lambda_0^2 - \gamma_0^2 = \lambda_0'^2+ \gamma_0'^2 \\
S_0 = 2 \lambda_0 \gamma_0 = 2 \lambda_0' \gamma_0'
\end{eqnarray}

Based on the previous identities, we can produce the following equalities:
\begin{eqnarray}\nonumber
\psi_0^2 = T_0^2 + S_0^2 \\
\phi_0^2 = T_0^2 - S_0^2
\end{eqnarray}

By multiplying the previous equalities, we have:  
\begin{equation}
x_0^2 = \psi_0^2 \phi_0^2 = T_0^4 - S_0^4
\end{equation} 

By adding and substracting them, we have 
\begin{equation}
2y_0^2 =  \frac{(\psi_0^2 + \phi_0^2)(\psi_0^2 - \phi_0^2)}{2} = 2T_0^2S_0^2
\end{equation}
By taking the sum of the squares of the two previous expressions we have: 
\begin{equation}
x_0^4 + 4y_0^4 = z_0^2 \text{ where } (x_0,2y_0)=1
\end{equation}

\subsection{$x^4 + 6x^2y^2 + y^4 = z^2$}

We shall now derive that the solvability of the equation $x^4 + 6x^2y^2 + y^4 = z^2$ (with the assumption $(x,y)=1$) is equivalent to the solvability of (\ref{resolvent1}); which in turn allows us to conclude the equivalence of solvability with $x^4 + 4y^4 = z^2$. 
\\
Suppose that we are given (\ref{resolvent1}) and it had solutions $(X_0,Y_0,X_0',Y_0')$, then we can use (\ref{resolvent1}) to construct the following 2-variable polynomial:
\begin{equation}
F(X,Y) = X^2 - Y^2 - (X'^2 + Y'^2) 
\end{equation}
Since $XY = X'Y'$ we can substitute : $Y^2 = (\frac{X'Y'}{X})^2$ resulting in 
\begin{equation}
G(X) = X^2  - (\frac{X'Y'}{X})^2 - (X'^2 + Y'^2)
\end{equation}
We can multiply by $X^2$ and we have:
\begin{eqnarray}
H(X) = X^2G(X) \\
H(X) = X^4  - (X'^2 + Y'^2)X^2 - (X'Y')^2
\end{eqnarray}

But, $H(X)$ is a quadratic polynomial over $X^2$, hence if it has integer roots (which it must have since $(X_0,Y_0,X_0',Y_0')$ are solutions of (8) ), its discriminant must have been a square in $\mathbb{Z}$. In other words:
\begin{eqnarray}
D^2 = (X'^2 + Y'^2)^2 + 4(X'Y')^2 \\
D^2 = X'^4 + 6X'^2Y'^2 + Y'^4
\end{eqnarray}
where if we set $X' = X_0'$ and $Y' = Y_0'$, this holds from some $D$. Hence we have shown that the solvability of (\ref{resolvent1}) implies the solvability of $x^4 + 6x^2y^2 + y^4 = z^2$.
\\
\\
The converse is also true. Suppose that indeed $(x_0,y_0,z_0)$ satisfied:
\begin{equation}
x^4 + 6x^2y^2 + y^4 = z^2 \hspace{20pt} (x,y)=1
\end{equation}

Then by subsitution  $t_0 = x_0^2 + y_0^2$, $s_0 = 2x_0y_0$, and using these new variables in the original equation we have $t_0^2 + s_0^2 = z_0^2$. By applying the pythagorean triplets, $t_0 = u_0^2 - v_0^2$ and $s_0 = 2u_0v_0$ which combined with our substitution yields the system of equations (8).

\section{Triviality of Solutions}

In the previous section we showed that $x^4 + 4y^4 = z^2$ and $x^4 + 6x^2y^2 + y^4 = z^2$ are equivalent to (\ref{resolvent1}). Now we will prove that (\ref{resolvent1}) has only trivial solutions over $\mathbb{Z}$.
\\
\\
Define:
\begin{equation}
\nu_p(a) = \alpha \text{ where } p^\alpha || a
\end{equation}
and
\begin{equation}
\nu(a) = \sum_{p \in \Pi} \nu_p(a)
\end{equation}

Let $(X_0,Y_0,X_0',Y_0')$ be a solution of (\ref{resolvent1}). By studying (\ref{resolvent1}) over modulo 2 and 3 respectively, we can conclude that $\nu(X_0Y_0) \geq 2$. Assume that $\nu(X_0Y_0)$ is minimal (i.e. there does not exist another $(X_1,Y_1,X_1',Y_1')$ such that it satisfies (\ref{resolvent1}) and $\nu(X_1Y_1) < \nu(X_0Y_0)$).
\\
\\
Since we have $X_0Y_0 = X_0'Y_0'$ where $(X_0,Y_0)=(X_0',Y_0')=1$, the following factorization is possible:
\begin{eqnarray*}
X_0 = (X_0,X_0')(X_0,Y_0') \\
Y_0 = (Y_0,X_0')(Y_0,Y_0') \\
X_0' = (X_0',X_0)(X_0',Y_0) \\
Y_0' = (Y_0',X_0')(Y_0',Y_0)
\end{eqnarray*}

Returning to (\ref{resolvent1}), re-write the first line as:

\begin{eqnarray}\nonumber
(X_0,X_0')^2(X_0,Y_0')^2 - (Y_0,X_0')^2(Y_0,Y_0')^2 \\
= (X_0,X_0')^2(X_0',Y_0)^2 + (Y_0',X_0)^2(Y_0,Y_0')^2
\end{eqnarray}

By re-arranging the terms to obtain products, we have:

\begin{eqnarray}\nonumber
[(X_0,X_0')^2 - (Y_0,Y_0')^2](X_0,Y_0')^2 \\
= [(X_0,X_0')^2 + (Y_0,Y_0')^2](X_0',Y_0)^2
\end{eqnarray}

Since we have $(X_0,Y_0')$ and $(X_0',Y_0)$ relatively prime, and $(X_0,X_0')^2 - (Y_0,Y_0')^2$ is relatively prime with $(X_0,X_0')^2 + (Y_0,Y_0')^2$, the equality could only hold if:

\begin{equation}
\begin{cases}
(X_0,Y_0')^2  = (X_0,X_0')^2 + (Y_0,Y_0')^2 \\
(X_0',Y_0)^2 = (X_0,X_0')^2 - (Y_0,Y_0')^2
\end{cases}
\end{equation}

But these are again primitive pythogorean triplets, hence the following representation must have been possible for some $s_0,t_0, s_0', t_0'$:

\begin{equation}
\begin{cases}
(X_0,X_0') = s_0^2 - t_0^2 \\
(Y_0,Y_0') = 2s_0t_0 \\
(X_0,X_0') = s_0'^2 + t_0'^2 \\
(Y_0,Y_0') = 2s_0't_0'
\end{cases}
\end{equation}

Which would in turn simply reduce to:

\begin{equation}
\begin{cases}
s_0^2 - t_0^2  = s_0'^2 + t_0'^2  \\
s_0t_0  = s_0't_0' \\
(s_0,t_0)  = (s_0',t_0')=1
\end{cases}
\end{equation}

Or be solution of (\ref{resolvent1}) where $\nu(s_0t_0) < \nu(X_0Y_0)$ a contradiction.

\section{Other Results and Discussion}

The experimental data that was collected has actually produced table 1 which shows that several quartics are equivalent to the same resolvent. However we would like to remark that there are equations, for example such as $x^4 + 2y^4 = z^2$ (whose triviality of solutions is shown in \cite{Pocklington1914} and \cite{Cohen2007}), which don’t belong to the below listed resolvents.

\begin{table}[H]
\centering
\scalebox{0.9}{
\begin{tabular}{|c | c|}
\hline
$
\begin{cases}
x^2 - y^2 = x'^2+ y'^2 \\
xy = x'y' \\
(x,y) = (x', y') = 1
\end{cases}
$
&  
$
\begin{array}{c}
x^4 - y^4 = z^2 \\
x^4 + 4y^4 = z^2 \\
x^4 + y^4 = 2z^2 \\
x^4 \pm 6x^2 y^2+y^4=z^2
\end{array}
$
\\
\hline
$
\begin{cases}
x^2 - 2y^2 = x'^2+ 2y'^2 \\
xy = x'y' \\
(x,y) = (x', y') = 1
\end{cases}
$ 
& 
$
\begin{array}{c}
x^4 + y^4 = z^2\\
x^4 - 4y^4 = z^2\\
x^4 - y^4 = 2z^2\\
x^4 \pm 12x^2 y^2+4y^4=z^2
\end{array}
$
\\
\hline
Both Resolvents & $x^4 + y^4 = z^4$ \\
\hline
\end{tabular}
}
\caption{Quartic equations and their resolvent}
\label{table1}
\end{table}

For algebraic extension $\mathbb{Z}[i]$, J. Cross using infinite descent over norms showed that the equation $x^4 + y^4 = z^2$ continues to have trivial solution \cite{Cross1993}. In this sense, it is plausible to believe that the technique developed in this paper for rational integers can also be successfully applied to algebraic extensions.

\bibliographystyle{plain}
\bibliography{new_bib.bib}

\end{document}